\documentclass[aps,%
pre,%
twocolumn,%
showpacs,%
floatfix,%
]{revtex4}

\usepackage{epsfig}

\begin{document}

\author{Jens Christian Claussen$^*$}

%\email{claussen@theo-physik.uni-kiel.de}
\affiliation{Institut f\"ur Theoretische Physik und Astrophysik, Universit\"at Kiel, Leibnizstra{\ss}e 15, 24098 Kiel,
Germany
%,  {\small \tt claussen@theo-physik.uni-kiel.de}
}

\title{Time-evolution of the Rule 150 cellular automaton activity
from a Fibonacci iteration
}
% described by a generalized hyper-Fibonacci series

\date{October 17, 2004}

\begin{abstract}
The total activity of the single-seeded 
cellular rule 150 automaton does not follow a one-step iteration
like other elementary cellular automata, 
but can be solved as a two-step vectorial, or string, iteration,
which can be viewed as a generalization of Fibonacci iteration
generating the time series from a sequence of vectors of 
increasing length.
This allows to compute the total activity time series
more efficiently than by simulating the whole spatio-temporal process,
or even by using the closed expression.
\end{abstract}
 %%% also 45.70.Ht, 02.30.Nw
\maketitle
 \vspace*{-6ex} \section{Introduction} \vspace*{-1.2ex}
\indent
Since the coining paper of Bak, Tang, and Wiesenfeld \cite{btw},
there has been considerable interest in the long-time   
behavior of cellular automata, especially for occurrence
of long range correlations, and correspondingly for
power spectra exhibiting a power law decay,
which have  become a paradigm for complex dynamical systems in general
\cite{jensen}.
% (see also  \cite{claussen2004}).
Cellular automata 
\cite{wolfram,wolfram_rmp}
are a widely studied class of discrete dynamical systems
showing emergence of complex spatio-temporal patterns
from a simple dynamical rule.
%
%
%\vspace*{-1.2ex} \section{ECA and sum signals} \vspace*{-1.2ex}
\paragraph*{ECA and sum signals} ---
A cellular automaton consists of an infinite
lattice of cells of e.g.\ two possible states (0,1), 
and a local deterministic update rule. At each discrete time step, 
a cell is updated $x_n^{t} \rightarrow x_n^{t+1}$ according to the
state within a local neighborhood.
For Conway's Game of Life \cite{conway} the
$3\times{}3$ Moore neighborhood on a two-dimensional lattice is used.
A simpler, yet complex, class are the 
elementary cellular automata (ECA)
\cite{wolfram,wolfram_rmp},
defined on a one-dimensional lattice, and the update rule
depends 
on the next-neighbor sites and its own state one time step
before:
\begin{eqnarray}
x_n^{t+1}=f(x_{n+1}^t, x_n^t, x_{n-1}^t)
\label{update90150}
\end{eqnarray}
where $f$ (the rule) is determined by $8$ bits being the output of the
possible
input bits $000$, $001$, ..., $111$;
this 8-bit number forms the ``rule'' number which enumerates the 
256  possible ECA rules.
The power spectra of rule 90 \cite{claussen2004} and some of the other ECA rules 
\cite{naglerclaussen04} exhibit a $1/f^\alpha$ decay.
Rule 90 and rule 150 can be expressed also as
\begin{eqnarray}
x_n^{t+1}=[x_{n-1}^t + r x_n^t + x_{n+1}^t] \mbox{\rm mod} 2
\end{eqnarray}
where $r=0$ defines rule 90, and $r=1$ rule 150, respectively.
In the context of catalytic processes, both can be interpreted 
as local self-limiting reaction processes
\cite{claussen2004,naglerclaussen04}.
As models for chemical turbulence,
similar dynamics with a continuous phase variable
have been discussed in 
\cite{twbcl}
and
\cite{oono},
including solitonic behavior, periodic, and turbulent
states. 
\paragraph*{Total activity.} ---
In the chemical picture, the total reaction rate at a given time $t$
corresponds to the total number of sites with $x_m^t=1$, 
described by the sum signal
\begin{eqnarray}  
X(t)=\sum_n x_n^t.
\end{eqnarray}
While for rule 150 there is no convenient solution
of $X(t)$ except a formal one
\footnote{
%\paragraph*{Closed solution from generating functions.~---} 
%
%The time evolution of the rule 150 ECA is given by \begin{eqnarray} x_n^{t+1}=(x^t_{n-1} + x^t_{n} +x^t_{n+1}) {\rm mod} 2. \end{eqnarray}
%
Closed solution from generating functions:
A straightforward expression 
for the time-evolution of $x^t_{n}$ 
by rule 150
for the initial condition of a single 1 at $n=0$ and $t=0$,
is defined by the 
coefficients of 
$P_t(x) := (1+x+x^2)$ \cite{wolfram_rmp},
or Gegenbauer functions,
%wolfram_rmp
%%%%%%% also known as Gegenbauer or ultraspherical functions,
and perform the mod 2 operation afterwards,
i.e.\ 
$P_t(x)=\sum_{n=0}^{2t} p_n  x^n$
and $x^t_{n}= p_{n-t}  {\rm mod} 2$.
From $P_t(x)$, the coefficients can be extracted
via the method of generating function,
giving
\centerline{$
%\begin{eqnarray}
x^t_{n}= \left[
\left.
\frac{\partial^{n+t}}{\partial x^{n+t}}
P_t(x)
\right|_{x=0}
\right] {\rm mod}  \; 2
%\end{eqnarray}
$}
Hovever in contrary to the case of rule 90, generating the
Sierpinski gasket, and where the cofficients
of $(1+x)^t$ simply are the binomial coefficients,
this does not lead to a closed expression for $X(t)$ rule 150.
},
and a fairly complicated expression 
(see Sec.~\ref{secexakt}),
it is computationally quite costly to
perform the full spatiotemporal dynamical simulation,
even if one is interested only in the time series.
This paper gives an iterative solution of $X(t)$
from a geometrical iteration and investigates the relationship to the 
Fibonacci iteration.
The block sums over 
%$0\ldots 2^N-1$ periods 
$0\leq t \leq 2^N-1$
can even be expressed 
directly via Fibonacci numbers.
\\
\indent
Throughout this paper the pure pattern
generated by a single 1
are considered on an infinite lattice.
%
%\clearpage \begin{widetext}
%
%
%
\vspace*{-1.2ex}
\section{Exact solution} \label{secexakt}
\vspace*{-1.2ex}
While in the Sierpinski (rule 90) case
$X(t)$ factorizes in a product of
$X_i(\sigma_i)$ for all ``time spins''
\cite{claussen2004},
for rule 150 it does not.
But again, a ``spin decomposition'' of time 
$t=\sum_{j=0}^{N-1}\sigma_j2^j$ with $\sigma_j\in\{0,1\}$
can be utilized as an efficient
coordinate system for the time axis.

Before turning to the geometric iteration,
it should be mentioned that a closed expression
in fact can be written down as follows.
As pointed out by Wolfram \cite{wolfram_rmp},
for rule 150 the ``correlation'' of the time spins comes into
play, i.e.\ $X(t)$ is exactly multiplicative for
blocks of spins of value 1  which are separated by one or more zeroes.
Then $X(t)=\prod_{n=1}^N  \chi(n)^{c_n}$,
where $c_n$ is the multiplicity of blocks of length $n$.
The series $\chi(n)$ should correctly read
 \footnote{The iteration  in
 \cite{wolfram_rmp} p.\ 614 gives wrong numbers for $\chi(n)$.}
to the iteration $\chi(n)=2\chi(n-1)-(-1)^n$ for $n\leq 1$
and $\chi(0)=1$.
Obviously $\chi(n)=X(2^N-1)$ holds.
The $\chi(n)$ on the other hand turn out as the
most decaying frequencies in the 
spectrum of rule 90
(see Fig.~2 in Ref.\ \cite{claussen2004}),
and can be expressed in closed form by
\begin{eqnarray}
\chi(n)=\left\lfloor \frac{2^{n+2}+1}{3} \right\rfloor,
\label{chiexakt}
\end{eqnarray}
where $\lfloor~\rfloor$ is the floor function.
%, ordownwards rounded integer.
%\\
%\clearpage
Defining $\sigma_{-1}:=0$ and $\sigma_{N}:=0$,
one can formalize the spin-block counting as
\begin{eqnarray}
X(t)\!=\!\prod_{n=1}^N \chi(n)^{^{
\displaystyle
\sum_{i=0}^{N-n}(1\!-\!\sigma_{i-1})(1\!-\!\sigma_{i+n})
\prod_{l=0}^{n-1}\sigma_{i+l}}}.
\end{eqnarray}
%  \clearpage
With our expression (\ref{chiexakt}),
this is a closed solution, and corresponds to
$X_{90}(t)=\prod_{j=0}^{N-1}2^{\sigma_j}$
in the Sierpinski case.
Due to the complicated time spin correlations, it
however looks quite unwieldy for analytical use,
and even is numerically unfeasible
 \footnote{Of order $o(t^3)$, optimized $o(t^2)$; the
 spatial simulation is also $o(t^2)$; 
 the iteration (\ref{zweimalvier}) is $o(t)$
 as the matrices are diagonal.}.
\vspace*{-1ex}
\section{Iterative solution by generalized hyper-Fibonacci series}
\vspace*{-1ex}
In contrast to the rule 90 (Sierpinski) case, for rule 150
the time evolution does not follow the same type 
of initiator-generator mechanism as it is well known 
for fractal sets.
However, it is possible to define a 
geometric
or measure-theoretic
\cite{georgeWgrossmann}
iteration based on the last and the last but one
iterate, see Fig.~\ref{replicationfig}.
This corresponds to a difference equation
with the r.h.s.\  depending
on the last two time steps, and in fact, for
the total activity within $2^N$ time steps we will
derive a difference equation later.

%\clearpage
 \begin{widetext}
\noindent
\begin{figure}[htbp]
\noindent
\epsfig{file=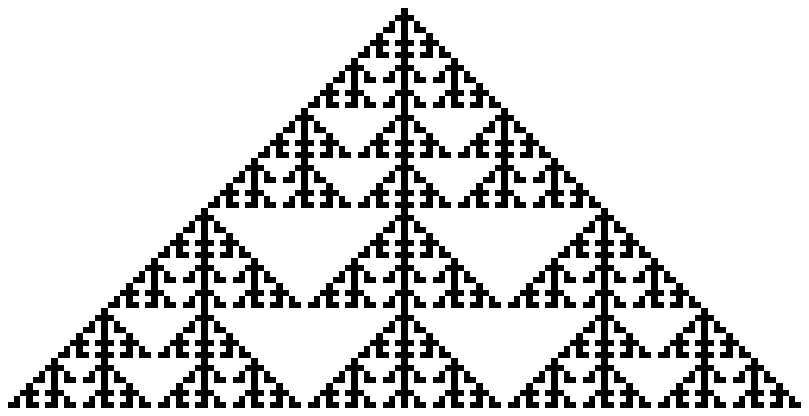,width=0.48\textwidth}
\hspace*{0.02\textwidth}
\epsfig{file=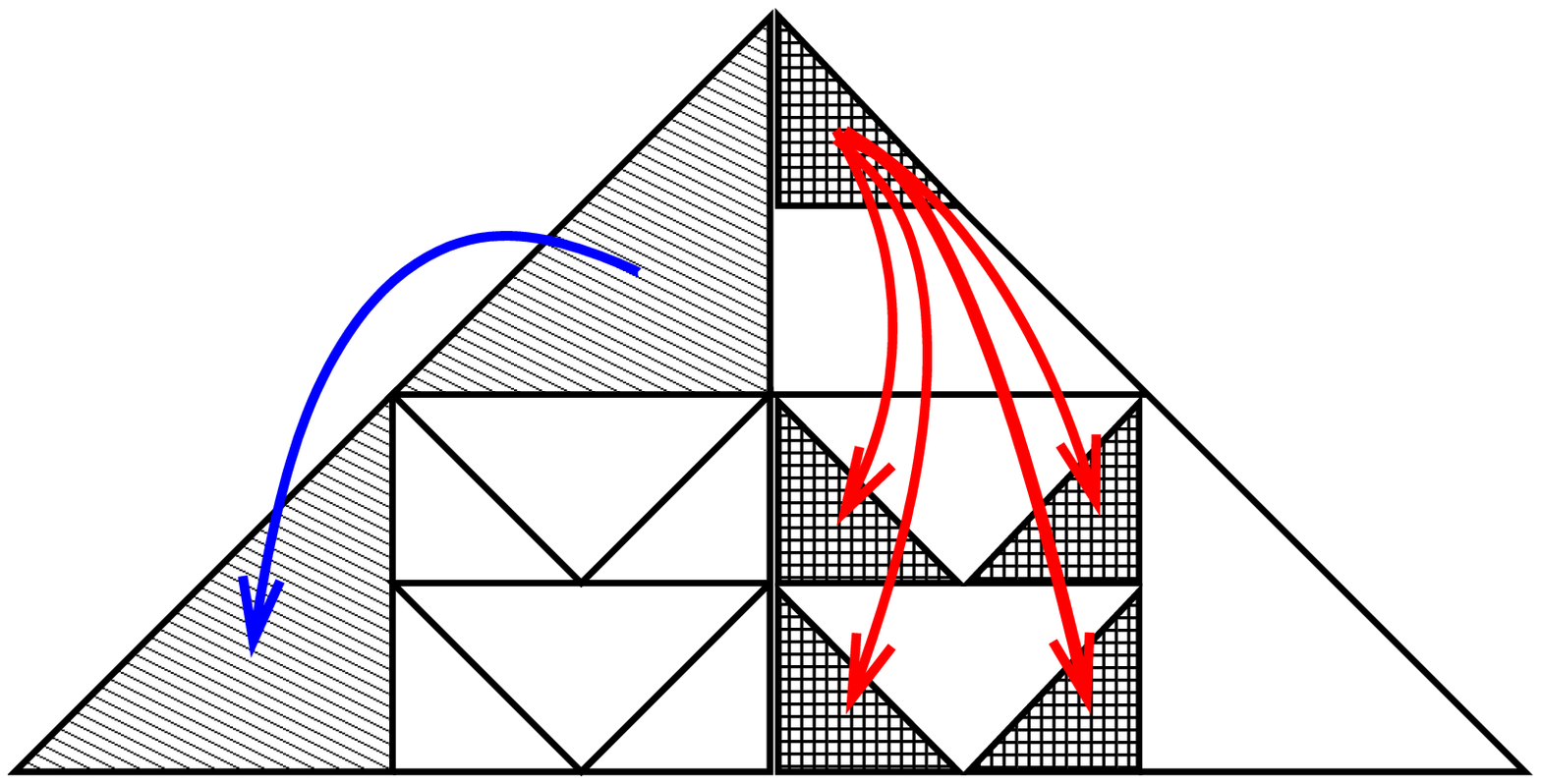,width=0.48\textwidth}
\caption{Left: Time-evolution of rule 150 for the first 64 time steps,
started with a single seed.
Right: Illustration of the replication rule. The whole system is
symmetric with respect to the vertical axis.
The whole triangle above is replicated once (left part of the triangle). 
The upper part is reproduced quadruplicate (right part of the triangle).
\label{replicationfig}
}
%\end{figure}
%
%\begin{figure}[htbp]
%
%\epsfig{file=tri150c.eps,width=0.48\textwidth}
%%%%\epsfig{file=rule150a.eps,width=\columnwidth}
%\caption{Illustration of the replication rule. The whole system is symmetric with respect to the vertical axis. The whole triangle above is replicated once (left).  The upper part is reproduced quadruplicate (right).}  \label{replicationfig}
\end{figure}

\noindent
According to the replication law
(see Fig. \ref{replicationfig}),
the time series of the total activity
$X(t)$ follows the two-step iteration 
%\begin{widetext}
\begin{eqnarray}
\label{x_it}
\nonumber
X(1,1,\sigma_{n-2} \cdots \sigma_0) 
&=& X(0,1,\sigma_{n-2} \cdots \sigma_0)
+ 2 X(0,0,\sigma_{n-2} \cdots \sigma_0)  \\
X(1,0,\sigma_{n-2} \cdots \sigma_0)&=& 
\hphantom{X(0,1,\sigma_{n-2} \cdots \sigma_0)}
+3 X(0,0,\sigma_{n-2} \cdots \sigma_0).
%\nonumber
\end{eqnarray}
%\end{widetext}
If we define
\begin{eqnarray}
Y(\bullet, ---) &:=& X(\bullet,0,---)
\nonumber
\\
Z(\bullet, ---) &:=& X(\bullet,1,---),
\end{eqnarray}
the iteration reads
%\begin{widetext}
\begin{eqnarray}   
\nonumber
Z(1,\sigma_{n-2} \cdots \sigma_0)
&=& Z(1,\sigma_{n-2} \cdots \sigma_0)
+ 2 Y(1,\sigma_{n-2} \cdots \sigma_0)
\\ 
Y(1,\sigma_{n-2} \cdots \sigma_0)&=&
\hphantom{Z(0,\sigma_{n-2} \cdots \sigma_0)}
+3 Y(1,\sigma_{n-2} \cdots \sigma_0),
%\nonumber
\end{eqnarray}
 \end{widetext}
or short
\begin{eqnarray}
\mathbf{Z}_n^0&=&3 \mathbf{Y}_{n-1}
\nonumber
\\
\mathbf{Z}_n^1
&=&
2 \mathbf{Y}_{n-1}
+ \mathbf{Z}_{n-1},
\end{eqnarray}  
and the concatenations 
% $Z_n=(Z_n^0,Z_n^1)$ and $Y_n=(Y_{n-1},Z_{n-1})$
\begin{eqnarray}
\mathbf{Y}_n&=&(\mathbf{Y}_{n-1},\mathbf{Z}_{n-1})^{\mbox{\scriptsize\tiny\rm T}}
\\
\mathbf{Z}_n&=&(\mathbf{Z}_n^0,\mathbf{Z}_n^1)^{\mbox{\scriptsize\tiny\rm T}}
\end{eqnarray}  
complete the iteration
\begin{eqnarray}
\label{y_it}
\left(
\mathbf{Y}_n
\!\!\!
\begin{array}{c}
\mbox{}
\\
\mbox{}
\end{array}
\right)
&=&
   \hphantom{\left(\begin{array}{cc}3\cdot \mbox{\large\bf 1} & \mbox{\large\bf 0}\\2\cdot \mbox{\large\bf 1} & \mbox{\large\bf 1}\end{array}\right) }
\left(
\begin{array}{c}
\mathbf{Y}_{n-1}\\
\mathbf{Z}_{n-1}\end{array}\right)
\\
\mbox{and} 
~~~~~~~
\left(
\mathbf{Z}_n
\!\!\!
\begin{array}{c}
\mbox{}
\\
\mbox{}
\end{array}
\right)
&=&
\left(
\begin{array}{cc}
3\cdot \mbox{\large\bf 1} & \mbox{\large\bf 0}\\
2\cdot \mbox{\large\bf 1} & \mbox{\large\bf 1}
\end{array}
\right)
\left(
\begin{array}{c}
\mathbf{Y}_{n-1}\\
\mathbf{Z}_{n-1}
\end{array}
\right)
%\end{eqnarray}
\label{z_it}
\end{eqnarray}
The initial vector is given by
$
\left(
\begin{array}{c}
\mathbf{Y}_{0}\\
\mathbf{Z}_{0}\end{array}\right)
=
\left(\begin{array}{c}1\\3\end{array}\right).$
Eqns.\ (\ref{y_it}-\ref{z_it}) can be collected together to the iteration
\begin{eqnarray}
\left(
\begin{array}{c}
\mathbf{Y}_n\\
\\
\mathbf{Z}_n
\end{array}
\!\!\!
\begin{array}{c}
\mbox{}
\\
\mbox{}
\\
\mbox{}
\\
\mbox{}
\end{array}
\right)
=
\left(
\begin{array}{cc}
 \mbox{\large\bf 1} & \mbox{\large\bf 0}\\
 \mbox{\large\bf 0} & \mbox{\large\bf 1}\\
3\cdot \mbox{\large\bf 1} & \mbox{\large\bf 0}\\
2\cdot \mbox{\large\bf 1} & \mbox{\large\bf 1}
\end{array}   
\right)
\left(
\begin{array}{c}
\mathbf{Y}_{n-1}\\
\mathbf{Z}_{n-1}
\end{array}   
\right)
\label{zweimalvier}
\end{eqnarray}
where the dimension of the vectors $\mathbf{Y}_n, \mathbf{Z}_n$ is  $2^n$,
growing in the same way as for  the Sierpinski iteration
(see Fig.~\ref{replicationfig90})
$
\left(
\mathbf{Z}_n
\!\!\!
\begin{array}{c}
\mbox{}
\\
\mbox{}
\end{array}
\right)
=\left(
\begin{array}{c}
 \mbox{\large\bf 1} \\ 2\cdot \mbox{\large\bf 1}
\end{array}
\right)
=\left( 
\begin{array}{c}
\mathbf{Y}_{n-1}\\
\mathbf{Z}_{n-1}\end{array}\right)$,
or the Thue-Morse iteration
$
\left(
\mathbf{Z}_n
\!\!\!
\begin{array}{c}
\mbox{}
\\
\mbox{}
\end{array}
\right)
=\left(
\begin{array}{c}
\mbox{\large\bf 1}\\
(-1)\cdot \mbox{\large\bf 1} 
\end{array}
\right)
\left( 
\begin{array}{c}
\mathbf{Y}_{n-1}\\
\mathbf{Z}_{n-1}\end{array}\right)$
(both use (\ref{y_it}) and start with $\mathbf{Z}_{0}=(1)$).

\begin{figure}[htbp]
\noindent
\epsfig{file=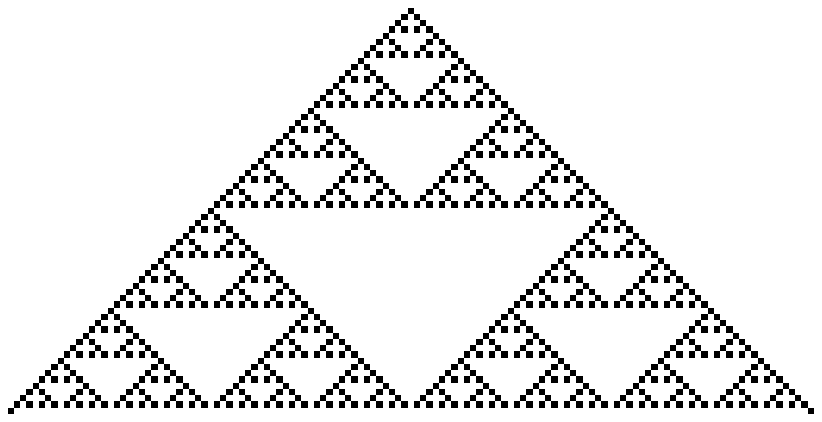,width=0.48\columnwidth}
\epsfig{file=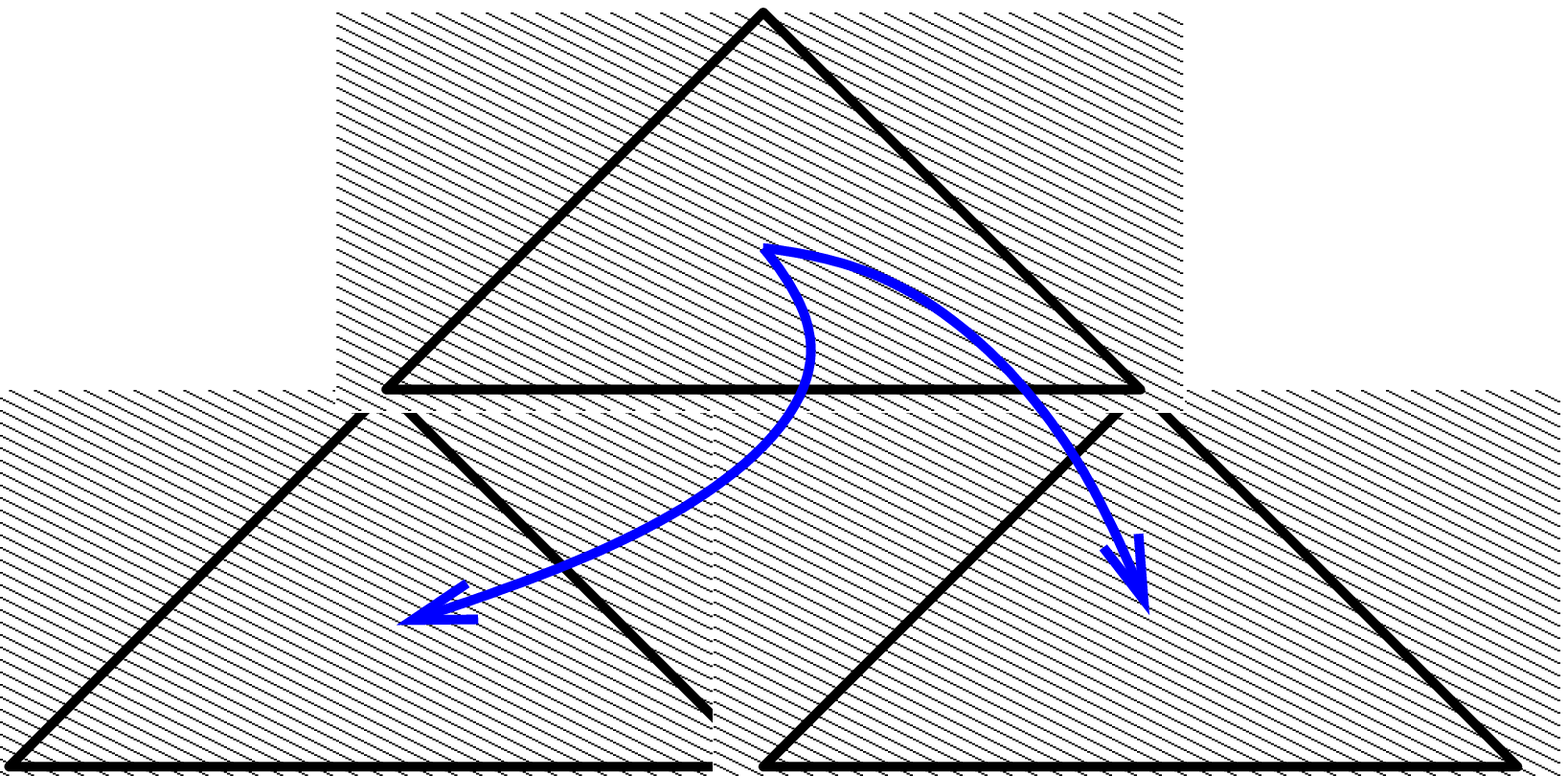,width=0.48\columnwidth}
%\caption{Time-evolution of rule 150 for the first 64 time steps, started with a single seed.}
%\epsfig{file=rule90_64s9.eps,width=0.48\columnwidth}
\caption{Illustration of the replication rule for rule 90.}
\label{replicationfig90}
\end{figure}

These iterations look formally similar to the
% (non-vectorial) 
Fibonacci or Lucas iteration
$\left(\begin{array}{c}F_{n}\\F_{n-1}\end{array}\right)
=
\left(\begin{array}{cc}
1&1\\1&0
\end{array}
\right)
\left(
\begin{array}{c}
F_{n-1}\\
F_{n-2}\end{array}\right)$.
Here 
%$ \left(\begin{array}{c} F_{1}\\F_{0}\end{array}\right) = \left(\begin{array}{c}  1\\0\end{array}\right) $
$ (F_{1},F_{0})^{\mbox{\tiny T}} =(1,0)^{\mbox{\tiny T}}  $
defines the Fibonacci series and
%$ \left( \begin{array}{c} F_{1}\\ F_{0}\end{array}\right) = \left( \begin{array}{c} 1\\2 \end{array}\right) $
$ (F_{1},F_{0})^{\mbox{\tiny T}} =(1,2)^{\mbox{\tiny T}}$
is the initial condition of the
Lucas series.
For the latter two, the length if the iterates 
is not growing.

Eqns.\ (\ref{x_it}) or (\ref{zweimalvier}), equivalently, together with the
initial condition $(1,3)^{\mbox{\tiny\rm T}}$
generate $X(t)$ iteratively for all t.
Formally this iteration is analogeous to the
Fibonacci iteration, but acts on vectors of growing length 
within an infinite-dimensional vector space indexed by
nonnegative integer values.

This type of series should be distinguished from the
($r$-th) hyper-Fibonacci series \cite{alexander_fibo},
where $f_l=2^{l-1}$ for 
$l=1,\cdots,r+1$ and $f_l=f_{l-1} + \cdots +
f_{l-r+1}$ for $l>r+1$.
On the other hand, the terminus {\sl generalized Fibonacci series} 
is widely used for the ordinary Fibonacci or Lucas iteration with
two arbitrary start values $f_0$ and $f_1$,
where $f_0=0$, $f_1=1$ defines the Fibonacci series
and $f_0=2$, $f_1=1$ defines the Lucas series;
and in fact both can be used as (nonorthogonal) basis vectors
of the linear space of generalized Fibonacci series.
The $r=1$ hyper-Fibonacci series corresponds to
a generalized Fibonacci series with  $f_0=1$, $f_1=2$.
Consequently, an iteration of the algebraic structure of 
Eq.\ (\ref{zweimalvier}) 
could be denoted as a generalized hyper-Fibonacci series.

\begin{figure}[htbp]
\epsfig{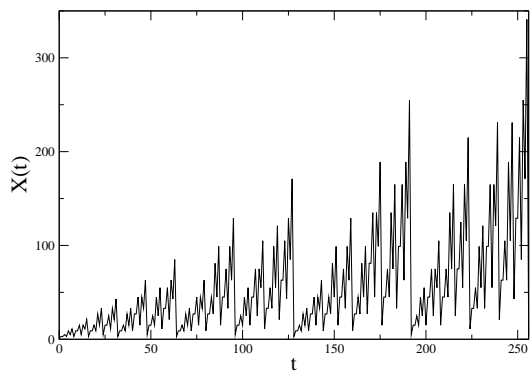}
\caption{\label{figx}
Rule 150: Plot of $X(t)$ for the first 256 time steps.
}
\end{figure}

\noindent 
Another observation is
the partial self-similarity relation
\begin{eqnarray}
X(\sigma_{n} \cdots\sigma_3,0,\sigma_1, \sigma_0) 
=
X(\sigma_1, \sigma_0)
\cdot
X(\sigma_{n} \cdots\sigma_3) 
\end{eqnarray}
(leading zeroes omitted in notation),
i.e.\ 
the sequence generated by every second four-block
$(\sigma_2=0)$ 
factorizes into the first block $(1,3,3,5)$
and the whole sequence itself.
A closed expression for 
$X(\sigma_{n} \cdots\sigma_3,1,\sigma_1, \sigma_0)$
is however not known yet.
\\
\indent
The first values of $X(t)=X(t_1+t_2)$ are listed in
Tab.~\ref{tabellex} (see Fig.~\ref{figx}).
\begin{table}[htbp]
%\vphantom{
\tiny\scriptsize
\hspace*{-1.0em}
\begin{tabular}{r|rrrrrrrrrrrrrrrr}
%$\downarrow\!\!t_1 \!\!\!$$\raisebox{1.8ex}{$t_2$}\!\!\!\!$ {{$\rightarrow$}}
&&&&&&&&&$t_2$\\
& 0 &16& 32&48 & 64&80
&96&112&128 &144& 160 &176& 192 &208& 224&240
%%%%%%%%%%%%%&256 &288 &320& 352 & 384 & 416 & 448 & 480 
\\
\hline
0 & 1 & 3 & 3 & 5 & 3 & 9 & 5 & 11 & 3 & 9 & 9 & 15 & 5 & 15 & 11 & 21  \\
1 & 3 & 9 & 9 & 15 & 9 & 27 & 15 & 33 & 9 & 27 & 27 & 45 & 15 & 45 & 33 & 63
 \\
2 & 3 & 9 & 9 & 15 & 9 & 27 & 15 & 33 & 9 & 27 & 27 & 45 & 15 & 45 & 33 & 63
 \\
3 & 5 & 15 & 15 & 25 & 15 & 45 & 25 & 55 & 15 & 45 & 45 & 75 & 25 & 75 & 55
& 105  \\
4 & 3 & 9 & 9 & 15 & 9 & 27 & 15 & 33 & 9 & 27 & 27 & 45 & 15 & 45 & 33 & 63
 \\
5 & 9 & 27 & 27 & 45 & 27 & 81 & 45 & 99 & 27 & 81 & 81 & 135 & 45 & 135 &
99 & 189  \\
6 & 5 & 15 & 15 & 25 & 15 & 45 & 25 & 55 & 15 & 45 & 45 & 75 & 25 & 75 & 55
& 105  \\
$\!\!\!t_1\!$ 
 7 & 11 & 33 & 33 & 55 & 33 & 99 & 55 & 121 & 33 & 99 & 99 & 165 & 55 & 165 &
121 & 231  \\
8 & 3 & 5 & 9 & 11 & 9 & 15 & 15 & 21 & 9 & 15 & 27 & 33 & 15 & 25 & 33 & 43
 \\
9 & 9 & 15 & 27 & 33 & 27 & 45 & 45 & 63 & 27 & 45 & 81 & 99 & 45 & 75 & 99
& 129  \\
10 & 9 & 15 & 27 & 33 & 27 & 45 & 45 & 63 & 27 & 45 & 81 & 99 & 45 & 75 & 99
& 129  \\
11 & 15 & 25 & 45 & 55 & 45 & 75 & 75 & 105 & 45 & 75 & 135 & 165 & 75 & 125
& 165 & 215  \\
12 & 5 & 11 & 15 & 21 & 15 & 33 & 25 & 43 & 15 & 33 & 45 & 63 & 25 & 55 & 55
& 85  \\
13 & 15 & 33 & 45 & 63 & 45 & 99 & 75 & 129 & 45 & 99 & 135 & 189 & 75 & 165
& 165 & 255  \\
14 & 11 & 21 & 33 & 43 & 33 & 63 & 55 & 85 & 33 & 63 & 99 & 129 & 55 & 105 &
121 & 171  \\
15 & 21 & 43 & 63 & 85 & 63 & 129 & 105 & 171 & 63 & 129 & 189 & 255 & 105 &
215 & 231 & 341  
\end{tabular}
%}%vhphantom
\normalsize
\caption{\label{tabellex}
The total activity $X(t)$ for the first 256 time steps.
}
\end{table}

\section{Iterative solution via a formal language replication rule}
One can proceed similar to the Sierpinski (rule 90) case, where $X(t)$
can be generated from an initiator $a=(1)$ by the 
iterative replication
$(a) \to (a,2\cdot a)$.
For the Thue-Morse sequence, the 
iteration reads 
$(a) \to (a,(-1)\cdot a)$,
see \cite{pikovskyPRE95,zaksPRL96,zaksPRL}.
However, for rule 150 again the Fibonacci-like structure comes into play,
i.e.\ each iterate depends on two preceders 
for which
one has to require ${\rm len}(a)={\rm len}(b)$.
Then the iteration is defined by
\begin{eqnarray}
(a,b) \to (a,b,3a,2a+b).
\end{eqnarray}
This is equivalent to Eq.~(\ref{zweimalvier}).

\section{Block-sums and the Fibonacci series}
Following the same geometrical argument as for the
row sums,
the sum $S_n=\sum_{i=0}^{2^n - 1}X(i)$ is given by the iteration
$S_n-S_{n-1} = S_{n-1} + 4 S_{n-2}$, or
\begin{eqnarray}
S_n = 2 S_{n-1} + 4 S_{n-2},
\end{eqnarray}
and the first elements of the series are listed in
Tab.~\ref{tabelles}.
\begin{table}[hbtp]
%\tiny\scriptsize\footnotesize\small
\begin{tabular}{rr|rr|rr}
$n$ & $S_n$ &$n$ & $S_n$ &$n$ & $S_n$ \\
0 & 1&6 & 1344&12 & 1544192\\
1 & 4&7 & 4352&13 & 4997120\\ 
2 & 12&8 & 14080&14 & 16171008\\
3 & 40&9 & 45568&15 & 52330496\\
4 & 128&10 & 147456&16 & 169345024\\
5 & 416&11 & 477184&17 & 548012032
%%% \\ &&&&18 & 1773404160
\end{tabular}
\normalsize  
\caption{\label{tabelles}
The block sums $S_n$ for the first 18 time steps.
}
\end{table}

The matrix of the iteration in time-delayed coordinates
\begin{eqnarray}
\left(
\begin{array}{c}
S_{n}\\
S_{n-1}\end{array}\right)
=
\left(
\begin{array}{cc}
2&4\\1&0
\end{array}
\right)
\left(
\begin{array}{c}
S_{n-1}\\
S_{n-2}\end{array}\right)
\end{eqnarray}
has the eigenvalues
$\alpha_{1,2}=1\pm\sqrt[]{5}$, indicating that it differs from the
Fibonacci iteration matrix by an additional expansion factor of 2, i.e.\
by a suitable transformation
\begin{eqnarray}
\left(
\begin{array}{c}
S_{n}\\
2 S_{n-1}\end{array}\right)
=
\left(
\begin{array}{cc}
1&1\\1&0
\end{array}
\right)
\left(
\begin{array}{c}
2 S_{n-1}\\
2^2 S_{n-2}\end{array}\right)
\end{eqnarray}
it relates to the Fibonacci numbers in usual convention,
\begin{eqnarray}
S_n &=& F_{n+2} \cdot 2^{n}
\\&=& \frac{1}{\sqrt{5}} \left[
\left(\frac{1+\sqrt{5}}{2}\right)^n -
\left(\frac{1-\sqrt{5}}{2}\right)^n \right],
\end{eqnarray}
the latter following from Binet's formula.
A ``blockwise normalization'', or detrended signal,
as used in \cite{claussen2004},
can be acheived by subtracting 
\begin{eqnarray}
N_n=(S_n-S_{n-1})/2^{n-1} = 2  F_{n-2} -  F_{n-2}
\end{eqnarray}
from $X(t)$ for all $2^{n-1} \leq t \leq 2^n - 1$.
Then
$
X(t) - N_{\lfloor \log_2 t \rfloor}
$
defines a signal with a mean vanishing 
within each time interval $[0,2^N-1]$.
---
The first values of this series are 
$N_0, N_1, \ldots = 1,3,4,7,11,18, \ldots$.

An immediate side result is that $S_n$ defines
the total volume of sites in $(n,t)$ space having the
value 1, and this volume scales for $n\to\infty$ 
with the largest eigenvalue. 
Thus, if the time-doubling iteration is interpreted as 
generation rule of the resulting self-similar 
fractal (rescaled to the unit interval), 
its Hausdorff-Besicovic dimension is 
given by $(1+\sqrt[]{5})/2$.
%
%\clearpage
%
\vspace*{-1.2ex}
\section{Conclusions}  
\vspace*{-1.2ex}
The self-similarity structure of the rule 150 elementary cellular automaton 
generated space-time fractal is 
qualitatively different from the Sierpinski triangle generated by rule 90.
While the iteration itself generalizes the concept of a 
Fibonacci iteration to vectors of growing dimension,
the blockwise sum exactly is given by the Fibonacci series
multiplyed by a scaling factor $2^n$. 
The iteration rule for the total activity derived here 
allows to compute the total activity without simulating the
spatial dynamics, thus considerably eases 
the numerical computation.
%\mbox{}
%\\[2ex]
\\[1ex]
\mbox{$^*$Electronic address: {\small\tt claussen@theo-physik.uni-kiel.de}}

\begin{widetext}
2000
{\sl
 Mathematics Subject Classification.
}
37B15, % Cellular automata (Topological dynamics)
11B39, % Fibonacci and Lucas numbers and polynomials and generalizations
28A80, % Fractals 
11B50, % Sequences (mod $m$)
11B85, % Automata sequences
68Q45, % Formal languages and automata 
68Q80. % Cellular automata (Theory of computing)
\\[0.5ex]
% \\\vfill
\end{widetext}

\end{document}